\documentclass[12pt,a4paper,draft]{article}
\usepackage[cp1251]{inputenc}
\usepackage{amsmath,amsthm,amssymb}
\usepackage{enumerate}
\usepackage[OT2,OT1]{fontenc}
\newtheorem{theorem}{Theorem}
\newtheorem{lemma}[theorem]{Lemma}
\newtheorem{corollary}[theorem]{Corollary}
\newtheorem{proposition}[theorem]{Proposition}
\newtheorem{remark}[theorem]{Remark}
\newcommand{\N}{{\mathbb N}}
\newcommand{\IS}{\mathcal {IS}_n}
\newcommand{\ISA}{(\mathcal {IS}_n, *)}
\newcommand{\T}{\mathcal {T}_n}
\renewcommand{\phi}{\varphi}
\renewcommand{\epsilon}{\varepsilon}
\newcommand{\im}{\operatorname{im}}
\newcommand{\Image}{\operatorname{Im}}
\newcommand{\rank}{\operatorname{rank}}
\newcommand{\dom}{\operatorname{dom}}

\newcommand{\Nil}{\operatorname{Nil}}
\newcommand{\Ord}{\operatorname{Ord}}
\newcommand{\ord}{\operatorname{ord}}
\newcommand{\Mon}{\operatorname{Mon}}
\newcommand{\type}{\operatorname{type}}
\begin{document}

\title{Isolated and nilpotent subsemigroups in the variants of $\IS$}
\author{G.Y. Tsyaputa}
\date{}
\maketitle

\begin{abstract}
All isolated, completely isolated, and nilpotent subsemigroups in
the semigroup $\IS$ of all injective partial transformations of an
$n$-element set, considered as a semigroup with a sandwich
multiplication are described.
\end{abstract}

\section{Introduction and main definitions}\label{s1}
In the monograph \cite{Ly} Ljapin proposed some constructions for
the semigroups, a certain modification of which is the following.
Let $S$ be a semigroup. For a fixed $a\in S$ define an operation
$*_a$ via $x*_a y=xay$, $x,y\in S$. Obviously, this multiplication
$*_a$ is associative, therefore the set $S$ with respect to this
operation is the semigroup which is called \textit{ the semigroup
} $S$ \textit{ with a sandwich operation } or the \textit{ variant
} of $S$ and is denoted by $(S,*_a)$.

The variants of the classical transformation semigroups,
$\mathcal{IS}_n$, $\mathcal{T}_n$, $\mathcal{PT}_n$, are
interesting examples of this construction. In \cite{Sy},
\cite{Ts}, and \cite{Ts3}, criteria of isomorphisms of these
variants are detailed, and some of their properties are described.

Recall, a subsemigroup $T$ of $S$ is called \textit{ isolated},
provided that $x^k\in T$ for some $k\in\N$ implies $x\in S$ for
all $x,y\in S$; $T$ is called \textit{ completely isolated},
provided that $xy\in T$ implies $x\in T$ or $y\in T$ for all
$x,y\in S$. Note that every completely isolated subsemigroup is
isolated, while the contrary does not hold true. Isolated and
completely isolated subsemigroups of the variants of $\T$ are
described in~\cite{MT}.

A semigroup $S$ with a zero $0$ is called \textit{nilpotent of the
nilpotency degree } $k\geq 1$, provided that $x_1 x_2 \cdots
x_k=0$ for all $x_1 , x_2 ,\dots , x_k \in S$ and there exist
$y_1, y_2,\dots ,y_{k-1} \in S$ such that $y_1 y_2 \cdots y_{k-1}
\neq 0$. Denote $n(S)=k$. It is known\cite{Arb}, that a finite
semigroup $S$ is nilpotent if and only if each element of $S$ is
nilpotent, that is, for each $x\in S$ there exist $k\in \N$ such
that $x^{k}=0$ and $x^{k-1}\neq 0$.

In this paper we describe completely isolated, isolated and
nilpotent subsemigroups of the variants of $\IS$, the inverse
symmetric semigroup of all partial injective transformations of
the set $N=\{1,2,\dots, n\}$.

In particular, in section \ref{s2} we study isolated and
completely isolated subsemigroups, and we produce the full
description of nilpotent subsemigroups of the variants of $\IS$
with respect to a natural zero (a nowhere defined map) in section
\ref{s3}. In the proofs we use the technique presented in
\cite{GM2}.

For an arbitrary $\beta\in\IS$ denote $\dom(\beta)$ and
$\im(\beta)$ the domain and image of $\beta$, respectively. The
value $|\dom(\beta)|=|\im(\beta)|$ is called the range of $\beta$
and is denoted by $\rank(\beta)$.

 It is proved in \cite{Ts} that semigroups $(\IS,*_\alpha)$ and $(\IS,*_\beta)$
are isomorphic if and only if $\pi\alpha=\beta\tau$ for some
permutations $\pi,\tau\in \mathcal{S}_n$. Note if $\alpha$ is a
permutation then $(\IS,*_\alpha)$ is isomorphic to $\IS$, and the
subsemigroup construction of the latter is studied in details in,
for example,\cite{GK}, \cite{GM1}. Therefore, without loss of
generality, we may assume that the sandwich element $\alpha$ is a
non-identity idempotent in $\IS$.

Now let us fix an idempotent $\alpha\in\IS$. Set
$\ISA=(\IS,*_{\alpha})$, $\dom(\alpha)=A$, $\rank(\alpha)=l$,
$l<n$, and let $\alpha$, as an element of $\IS$, be defined on $A$
identically.

Fix some $z\in N\setminus A$ and for arbitrary pairwise different
elements $x_1,\dots,x_k\in A$ define the following partial
permutations:
\begin{eqnarray}\nonumber (x_1,\dots ,x_k) (x)=
 \left\{
 \begin{array}{lccr} x_{i+1}, \quad \text{ if } x=x_i,\quad i<k,
 \\
 x_1, \quad \text{ if } x=x_k,
 \\
 x, \quad \text{ if } x\in N\setminus\{x_1,\dots,x_k\}
 \end{array}
\right.
\end{eqnarray}
and
\begin{eqnarray}\nonumber [x_1,\dots ,x_k] (x)=
 \left\{
 \begin{array}{lccr} x_{i+1}, \quad \text{ if } x=x_i,\quad i<k,
 \\
 z, \quad \text{ if } x=x_k,
 \\
 x, \quad \text{ if } x\in N\setminus\{x_1,\dots,x_k, z\}.
 \end{array}
\right.
\end{eqnarray}
For any $\beta\in\ISA$ we denote $\beta^{*s}=\underbrace{\beta
*\beta*\dots *\beta}_s$.

\begin{proposition}[\cite{Ts}]
An element $\epsilon\in\IS$ is an idempotent in $\ISA$ if and only
if $\epsilon$ is an idempotent in $\IS$ and
$\dom(\epsilon)\subset A$.
\end{proposition}
\begin{remark}\label{id} For each $x\in A$ denote $\epsilon_x$ the idempotent of
$\IS$ such that $\dom(\epsilon_x)=A\setminus\{x\}$. Then any
idempotent $\epsilon\neq\alpha$ of $\ISA$ may be factorized as
\begin{displaymath}
\epsilon=\prod_{x\in
A\setminus\dom(\epsilon)}\epsilon_x=\underset{x\in
A\setminus\dom(\epsilon)}{{\prod}^*} \negthickspace\; \epsilon_x.
\end{displaymath}
\end{remark}

\section{Isolated and completely isolated subsemigroups in $\ISA$}\label{s2}
Let $S$ be a semigroup. For the idempotent $e$ in $S$ define
\begin{displaymath}
 \sqrt{e}=\{x\in S: x^m=e\text{ for some } m>0\}.
\end{displaymath}
\begin{proposition}[\cite{MT}]\label{sqrt}
If $\sqrt{e}$ is a subsemigroup of $S$ then $\sqrt{e}$ is a
minimal with respect to inclusion isolated subsemigroup containing
$e$.
\end{proposition}
Denote $\mathcal{C}_A$ the set of those elements from $\IS$, which
are one-to-one maps on $A$ and are arbitrarily defined on
$N\setminus A$. That is,
\begin{displaymath} \mathcal{C}_A
=\{\beta\in\mathcal{IS}_n \;|\; \beta(A)=A \}.
\end{displaymath}
\begin{theorem}\label{cisol} The only completely isolated subsemigroups of $\ISA$ are $\mathcal{C}_A$,
$\ISA\setminus\mathcal{C}_A$ and $\ISA$.
\end{theorem}
\begin{proof}
It is clear that $\mathcal{C}_A$ is a subsemigroup of $\ISA$. Let
$\beta *\gamma \in \mathcal{C}_A$ for some $\beta$, $\gamma$ from
$\ISA$. Then $\dom(\beta *\gamma)\supset A$, therefore $\beta
(x)\in A$ for all $x\in A$. This implies $\beta\in\mathcal{C}_A$.
Hence $\mathcal{C}_A$ is a completely isolated subsemigroup.
$\ISA\setminus\mathcal{C}_A$ is also completely isolated as a
complement to the completely isolated subsemigroup. Obviously,
$\ISA$ is completely isolated.

Conversely, let $T\subset\ISA$ be a completely isolated
subsemigroup. Assume that $T\cap\mathcal{C}_A \neq\varnothing$,
$\beta\in \mathcal{C}_A\cap T$. Then $\beta^{*s} =\alpha$ for some
$s\in\N$ hence $\alpha\in T$. However, for any $\gamma$ from
$\mathcal{C}_A$ we have $\gamma^{*t}=\alpha$ for some $t\in\N$,
hence $\gamma\in T$. Therefore
\begin{equation}\label{riv1}
T\cap\mathcal{C}_A =\varnothing \quad\text{or}\quad
\mathcal{C}_A\subset T.
\end{equation}
Let $T\cap(\ISA\setminus\mathcal{C}_A)\neq\varnothing$ and $\beta$
belong to the intersection. Then $T$ contains an idempotent
$\epsilon$ as a power of $\beta$ and $\rank (\epsilon)<l$. From
remark \ref{id} it follows that since $T$ is completely isolated
$\epsilon_x\in T$ for some $x\in A$. But for any $y\in A$ we have
$\epsilon_x=(x,y) *\epsilon_y* (x,y)$. So if $\mathcal{C}_A\subset
T$ then $T$ contains $\epsilon_y$, as it is a semigroup. Otherwise
$T$ contains $\epsilon_y$ as it is completely isolated.
Consequently, $T$ contains all idempotents of ranks $<l$, and
since some power of each element $\gamma\in
\ISA\setminus\mathcal{C}_A$ equals such idempotent, we get
\begin{equation}
\label{riv2} (\ISA\setminus\mathcal{C}_A)\subset T \quad
\text{or}\quad (\ISA\setminus\mathcal{C}_A)\cap T=\varnothing.
\end{equation}
Now the statement of the theorem follows from (\ref{riv1}) and
(\ref{riv2}).
\end{proof}
For every $x\in A$ denote
\begin{displaymath}
\mathcal{G}(x)=\{\lambda\in\IS\; : \lambda(x)\notin A \quad\text{
and }\quad \lambda(y)\in A\setminus\{x\}\quad \text{ for all }
y\in A\setminus\{x\}\}.
\end{displaymath}
Observe, if $|A|=1$ then
$\mathcal{G}(x)=\ISA\setminus\mathcal{C}_A$, $x\in A$.
\begin{lemma}\label{lemisol}
$\mathcal{G}(x)=\sqrt{\epsilon_x}$.
\end{lemma}
\begin{proof}
It is clear that $\mathcal{G}(x)$ is a subsemigroup of $\ISA$ and
$\mathcal{G}(x)\subset\sqrt{\epsilon_x}$. Further
$\epsilon_x\in\mathcal{G}(x)$ by definition of $\mathcal{G}(x)$.

Let $\lambda^{*k}\in \mathcal{G}(x)$ for some $k\in\N$. If
$\lambda(x)\in A$ then whereas $\lambda^{*k}$ is not one-to-one
map on $A$ there exists $y\in A\setminus\{x\}$ such that
$\lambda(y)\in N\setminus A$. Then $\lambda^{*k}(y)$ is not
defined for $k\geq 2$, which contradicts the definition of
$\mathcal{G}(x)$. Now let $\lambda(x)\in N\setminus A$. Clearly,
for $\lambda^{*k}\in\mathcal{G}(x)$ it is necessary $\lambda(y)\in
A\setminus\{x\}$ for all $y\in A\setminus\{x\}$. Hence,
$\lambda\in\mathcal{G}(x)$ and $\mathcal{G}(x)$ is isolated
subsemigroup. Finally, proposition \ref{sqrt} completes the proof.
\end{proof}
\begin{theorem}
\begin{enumerate}[(i)]
\item If $\rank(\alpha)\geq 2$ then the only isolated subsemigroups
of $\ISA$ are $\ISA$, $\mathcal{C}_A$, $\ISA\setminus
\mathcal{C}_A$ and $\mathcal{G}(x)$, $x\in A$.
\item If $\rank(\alpha)=1$ then the only isolated subsemigroups of
$\ISA$ are $\mathcal{C}_A$, $\ISA\setminus\mathcal{C}_A$, $\ISA$,
in particular all of them are completely isolated.
\end{enumerate}
\end{theorem}
\begin{proof}
By theorem \ref{cisol} and lemma \ref{lemisol} all listed
subsemigroups are isolated.

Let $T$ be isolated subsemigroup of $\ISA$. As in the proof of
theorem \ref{cisol} it can be shown that (\ref{riv1}) holds.

Assume that $\mathcal{C}_A\subset T$. If $T\neq\mathcal{C}_A$ then
there exists $\beta\in T\setminus\mathcal{C}_A$ and $T$ contains
some idempotent $\epsilon$ of rank $<l$. Denote
$A\setminus\im(\epsilon)=\{a_1,\dots,a_k\}$ and consider elements
$\lambda=[a_1,\dots,a_k]$ and $\mu=[a_k,\dots,a_1]$. From
$\lambda^{*k}=\mu^{*k}=\epsilon$ it follows that $\lambda,\mu \in
T$ and $\lambda *\mu=\epsilon_{a_k}$ belongs to $T$. However as
$\mathcal{C}_A\subset T$, $T$ contains all $\epsilon_{x}$, $x\in
A$, and hence it contains all idempotents of ranks $\leq l-1$.
Therefore in this case $T=\ISA$.

Now assume that $\mathcal{C}_A\cap T=\varnothing$ and $T$ contains
an element $\beta$ and some idempotent $\epsilon$ of rank $<l$ as
a power of $\beta$. Consider the cases.
\begin{enumerate}[1)]
\item The case of $l=1$.
\\
Obviously, the power of any element from $\ISA\setminus
\mathcal{C}_A$ is a nowhere defined map, so $T=\ISA\setminus
\mathcal{C}_A$.
\item The case of $l\geq 2$.

If $\rank(\epsilon)\leq l-2$ then by the analogous arguments
produced above it can be shown that $T$ contains at least two
different idempotents $\epsilon_x$ and $\epsilon_y$, $x,y\in A$,
$x\neq y$. We show that $T$ contains all idempotents of rank
$l-1$, and hence, $T=\ISA\setminus\mathcal{C}_A$. Indeed, let
$z\in A$ and $z\neq x,y$. Consider element $\mu=(x,z) *[y]$.
$\mu^{*2}=\epsilon_y$, hence $\mu\in T$. Then
$(\epsilon_{x}*\mu)^{*2}=[x]*[y]*[z] \: \in T$. Let
$\sigma=[x,y,z]$, $\tau=[z,y,x]$. We have
$\sigma^{*3}=\tau^{*3}=[x]*[y]*[z]$, hence $\sigma, \tau \in T$
and finally $\sigma *\tau =\epsilon_{z} \: \in T$.

Now let $\rank(\epsilon)=l-1$. If $\epsilon$ is not the only
idempotent in $T$ then by the above reasoning
$T=\ISA\setminus\mathcal{C}_A$. If $\epsilon$ is the only
idempotent in $T$ then $T=\mathcal{G}(x)$.
\end{enumerate}
The theorem is proved.
\end{proof}

\section{Nilpotent subsemigroups in $\ISA$}\label{s3}

For every positive integer $k$ denote $\Nil_k$ the set of all
nilpotent subsemigroups of $\ISA$ of nilpotency degree $\leq{k}$.
The set $\Nil_k$ is partially ordered with respect to inclusions
in a natural way. Set $M=\overline{A}^{(1)}\cup
A\cup\overline{A}^{(2)}$, where $\overline{A}^{(1)}=N\setminus A$,
and $\overline{A}^{(2)}$ is a disjoint copy of
$\overline{A}^{(1)}$. For every $x\in \overline{A}^{(1)}$ denote
$x'$ the corresponding element from $\overline{A}^{(2)}$.

Denote $\Ord_k (M)$ the ordered set of all strict partial orders,
$\Lambda$, on $M$ which satisfy the following two conditions:
\begin{enumerate}[(1)]
\item\label{t6.1.1} the cardinalities of chains of $\Lambda$ are bounded by $k$,
\item\label{t6.1.2} $\overline{A}^{(1)}\subseteq\min(\Lambda)$,
$\overline{A}^{(2)}\subseteq\max(\Lambda)$, where $\min(\Lambda)$
and $\max(\Lambda)$ mean the sets of all minimal and maximal
elements of order $\Lambda$ respectively.
\end{enumerate}
If $k\leq m$ then we have natural inclusions $\Nil_k
\hookrightarrow\Nil_m$ and $\Ord_k (M) \hookrightarrow\Ord_m (M)$,
which preserve the partial order. Therefore we can consider the
ordered sets
\begin{displaymath}
\Nil =\bigcup_{k}\Nil_k \quad \text{ and }\quad \Ord(M)
=\bigcup_{k}\Ord_k(M).
\end{displaymath}
For every partial order $\Lambda\in\Ord(M)$ consider the set
\begin{equation}\label{mon}
\Mon(\Lambda)=\{\beta\in\IS : \beta(x)\neq x \text{ and } (x,
\beta(x))\in\Lambda \text{ for all } x\in\dom(\beta)\}
\end{equation}
 and for every subsemigroup $S\in\Nil$ the relation
\begin{align}\label{ord}\nonumber
 \Lambda_S &=\{(x,y): x\in A \text{ and there exists } \beta\in S \text{ such that }
 \beta(x)=y\} \\\nonumber &\cup\{(x,y): y\in A \text{ and there exists } \beta\in S \text{ such that }
 \beta(x)=y\} \\&\cup\{(x,y): x\in\overline{A}^{(1)}, y\in\overline{A}^{(2)} \text{
 and there exists } \beta\in S \text{ such that }
 \beta(x)=y\}.
\end{align}

Let $\beta\in\IS$ and $\rank(\beta)=k$, $k\leq n$. Write $\beta$
as
\begin{equation}\label{beta}
\beta=\left(
\begin{array}{cccccccccccccccccccccccccccccccccccccccccc}
x_{1 1} &\dots& x_{1 i_1}&x_{2 1}&\dots &x_{2 i_2}& x_{3 1} & \dots & x_{3 i_3}&x_{4 1}&\dots & x_{4 i_4}\\
y_{1 1} &\dots& y_{1 i_1}&y_{2 1}&\dots &y_{2 i_2}& y_{3 1} &
\dots & y_{3 i_3}&y_{4 1}&\dots & y_{4 i_4}
\end{array}
\right),
\end{equation}
that is, $\dom(\beta)=\{x_{1 1},\dots ,x_{4 i_4}\}$ and
$\beta(x_{i j})=y_{i j}$ for all $i,j$, moreover
\begin{eqnarray*}
&&\{x_{1 1}, \dots, x_{1 i_1}\}\subset A,\qquad\qquad
\{y_{1 1}, \dots, y_{1 i_1}\}\subset A,\\
&&\{x_{2 1},\dots, x_{2 i_2}\}\subset A,\qquad\qquad
\{y_{2 1},\dots, y_{2 i_2}\}\subset N\setminus A,\\
&&\{x_{3 1},\dots, x_{3 i_3}\}\subset N\setminus A,\qquad
\{y_{3 1},\dots, y_{3 i_3}\}\subset A,\\
&&\{x_{4 1},\dots, x_{4 i_4}\}\subset N\setminus A,\qquad \{y_{4
1}, \dots, y_{4 i_4}\}\subset N\setminus A.
\end{eqnarray*}
Consider a map $f: \ISA \rightarrow \mathcal{IS}(M)$ defined in
the following way:

if $\beta$ is given by (\ref{beta}) then $\dom(f(\beta))=\{x_{1
1},\dots ,x_{1 i_1}, x_{2 1},\dots ,x_{4 i_4}\}$ and
\begin{equation*}\label{f(beta)}
f(\beta)=\left(
\begin{array}{cccccccccccccccccccccccccccccccccccccccccc}
x_{1 1} &\dots& x_{1 i_1}&x_{2 1}&\dots &x_{2 i_2}& x_{3 1} & \dots & x_{3 i_3}&x_{4 1}&\dots & x_{4 i_4}\\
y_{1 1} &\dots& y_{1 i_1}&y'_{2 1}&\dots &y'_{2 i_2}& y_{3 1} &
\dots & y_{3 i_3}&y'_{4 1}&\dots & y'_{4 i_4}
\end{array}
\right).
\end{equation*}
\begin{proposition}\label{mono}
The above map $f: \ISA \rightarrow \mathcal{IS}(M)$ is a
monomorphism, besides,
\begin{displaymath}\Image(f)=\{\gamma\in \mathcal{IS}(M) :
\dom(\gamma)\subset \overline{A}^{(1)}\cup A, \im(\gamma)\subset
A\cup\overline{A}^{(2)}\}.
\end{displaymath}
\end{proposition}
\begin{proof}Clear from the definition.
\end{proof}
Hence, every nilpotent subsemigroup of $\ISA$ is mapped by $f$ to
the corresponding nilpotent subsemigroup of $\mathcal{IS}(M)$.
This allows one to apply the results from \cite{GM2} for
$\mathcal{IS}(M)$ to the semigroup $\ISA$. In particular, by a
word for word repetition of corresponding proofs from \cite{GM2}
one may prove propositions \ref{tv7}-\ref{tv10}:
\begin{proposition}\label{tv7}
\begin{enumerate}[(i)]
\item\label{tv7.1} For each $k\geq 1$ the map
$\Lambda\mapsto\Mon(\Lambda)$ is a homomorphism from the poset
$\Ord_k (M)$ to the poset $\Nil_k$.
 \item\label{tv7.2} For every $k\geq 1$ the map
$S\mapsto\Lambda_S$ is a homomorphism from the poset $\Nil_k$ to
the poset $\Ord_k (M)$.
\end{enumerate}
\end{proposition}

\begin{proposition}\label{tv9}
Let $n(S)=k$. Then $\Lambda_S \in
\Ord_k(M)\setminus\Ord_{k-1}(M)$.
\end{proposition}

\begin{proposition}\label{tv10}
Let $S\in\Nil$ and $\Lambda\in\Ord(M)$. Then
\begin{enumerate}[(i)]
\item\label{tv10.1} $\Mon(\Lambda_S)\supset S$,
$\Lambda_{\Mon(\Lambda)}\subset\Lambda$;
\item\label{tv10.2} $\Mon(\Lambda_{\Mon(\Lambda)})=\Mon(\Lambda)$;
\item\label{tv10.3} $\Lambda_{\Mon(\Lambda_S)}=\Lambda_S$.
\end{enumerate}
\end{proposition}
From proposition \ref{tv9} we derive
\begin{corollary}
The nilpotency degree, $n(S)$, of any nilpotent subsemigroup
$S\subset\ISA$ does not exceed $\rank(\alpha)+2$. Nilpotent
subsemigroup $S$ has nilpotency degree $\rank(\alpha)+2$ if and
only if $|N\setminus\im(\alpha)|\geq 1$,
$\min(\Lambda_S)=\overline{A}^{(1)}$,
$\max(\Lambda_S)=\overline{A}^{(2)}$ and the restriction of
$\Lambda_S$ on $\im(\alpha)$ is a linear order.
\end{corollary}

By an \textit{ ordered } $A$\textit{-partition } of $M$ into $k$
non-empty blocks we mean the partition $M=M_1\cup M_2\cup\dots\cup
M_k$ where $M_1\supset\overline{A}^{(1)}$,
$M_k\supset\overline{A}^{(2)}$, and the order of the blocks is
also taken into account. With every ordered $A$-partition of
$M=M_1\cup M_2\cup\dots\cup M_k$ we associate the set
\begin{equation}\label{maxord}
\ord(M_1,\dots,M_k)=\bigcup_{1\leq i<j\leq k} M_i \times M_j
\subset M\times M.
\end{equation}
\begin{lemma}\label{max}
Let $k\leq |M|$ be fixed. Then
\begin{enumerate}[(i)]
\item\label{max1} for every ordered $A$-partition $M=M_1\cup M_2\cup\dots\cup
M_k$ the set $\ord(M_1,\dots,M_k)$ is a maximal element in $\Ord_k
(M)$;
\item\label{max2} different ordered $A$-partitions of $M$ correspond to
different elements in $\Ord_k (M)$;
\item\label{max3} each maximal element in
$\Ord_k (M)$ has the form $\ord(M_1,\dots,M_k)$ for some ordered
$A$-partition of $M=M_1\cup M_2\cup\dots\cup M_k$.
\end{enumerate}
\end{lemma}
\begin{proof}Statements \eqref{max1} and \eqref{max2} are proved analogously to
lemma 7 from \cite{GM2}.
\\ \eqref{max3} Let the order $\Lambda\in\Ord_k (M)$ be fixed.
Denote by $M_1$ the set of all minimal elements of the relation
$\Lambda$. By definition $\overline{A}^{(1)}\subset M_1$,
therefore $M_1\neq \varnothing$. For every increasing (with
respect to $\Lambda$) chain $x_1<\dots <x_m$ of elements in
$M\setminus M_1$ there exists $x_0\in M$ such that $x_0<x_1$.
Hence, the cardinality of every increasing chain in $M\setminus
M_1$ is bounded by $k-1$. Denote by $M_2$ the set of all minimal
elements in $M\setminus M_1$, by $M_3$ the set of all minimal
elements in $M\setminus (M_1\cup M_2)$ and so on. In $k$ steps we
get the partition $M=M_1\cup M_2\cup\dots\cup M_k$, for which
$\overline{A}^{(1)}\subset M_1$. Observe that the elements of
$\overline{A}^{(2)}$ are maximal by the definition of $\Lambda$.
Hence set
\begin{eqnarray*}
&&\widetilde{M_i}=M_i\setminus \overline{A}^{(2)}, \quad 1\leq
i\leq
k-1,\\
 &&\widetilde{M_k}=M_k\cup \overline{A}^{(2)}
\end{eqnarray*}
and get the new partition $M=\widetilde{M_1}\cup
\widetilde{M_2}\cup\dots\cup \widetilde{M_k}$. Clearly,
$\Lambda\subset\ord(\widetilde{M_1},\dots, \widetilde{M_k})$ and
from the maximality of $\Lambda$ we have $\Lambda
=\ord(\widetilde{M_1},\dots, \widetilde{M_k})$. The lemma is
proved.
\end{proof}
For every ordered $A$-partition of $M=M_1\cup M_2\cup\dots\cup
M_k$ denote
\begin{multline}\label{maxnil}
T(M_1,\dots,M_k)=\{\beta\in\IS : x\in M_i \text{ and } \beta(x)\in
M_j \text{ imply } i<j \\ \text{ for all } x\in\dom(\beta)\}.
\end{multline}
\begin{lemma}\label{rem2}
\begin{displaymath}
T(M_1,\dots,M_k)=\Mon(\ord(M_1,\dots,M_k)).
\end{displaymath}
\end{lemma}
\begin{proof}
The assertion follows from definitions (\ref{mon}), (\ref{ord}),
(\ref{maxord}) and (\ref{maxnil}).
\end{proof}
\begin{theorem}\label{teor14}
\begin{enumerate}[(i)]
\item\label{14.1} For every ordered $A$-partition of $M=M_1\cup
M_2\cup\dots\cup M_k$ the semigroup $T(M_1,\dots,M_k)$ is maximal
in $\Nil_k$.
\item\label{14.2} Different ordered $A$-partitions of $M$ correspond to different subsemigroups in
$\Nil_k$.
\item\label{14.3} Every maximal subsemigroup in $\Nil_k$ has the form $T(M_1,\dots , M_k)$
for some ordered $A$-partition of $M=M_1\cup M_2\cup\dots\cup
M_k$.
\end{enumerate}
\end{theorem}
\begin{proof}
From lemma \ref{rem2} it follows that the set $T(M_1,\dots,M_k)$
is a subsemigroup in $\Nil_k$. We show that
\begin{equation}\label{rivnist}
\Lambda_{T(M_1,\dots,M_k)}=\ord(M_1,\dots,M_k).
\end{equation}
From proposition \ref{tv10}\eqref{tv10.1} we have
$\Lambda_{T(M_1,\dots,M_k)}\subset\ord(M_1,\dots,M_k)$.
\\ We prove the contrary inclusion. Let $(x,y)\in
\ord(M_1,\dots,M_k)$. Take $\beta\in\IS$ such that
$\rank(\beta)=1$, $\dom(\beta)=\{x\}$, $\im(\beta)=\{y\}$. By
definition $\beta\in\Mon(\ord(M_1,\dots,M_k))$, and using
(\ref{ord}) we have $(x,y)\in\Lambda_{\Mon(\ord(M_1,\dots,M_k))}$.
Hence $\Lambda_{T(M_1,\dots,M_k)}\supset\ord(M_1,\dots,M_k)$ and
the equality (\ref{rivnist}) is proved.

Now let $S\in\Nil_k$ be such that $S\supset T(M_1,\dots,M_k)$.
According to lemma \ref{max} the order $\ord(M_1,\dots,M_k)$ is a
maximal element in $\Ord_k (M)$, therefore by proposition
\ref{tv7}\eqref{tv7.2} we get $\Lambda_S=\ord(M_1,\dots,M_k)$.
Finally from proposition \ref{tv10}\eqref{tv10.1}, lemma
\ref{rem2}, and equality (\ref{rivnist}) it follows
\begin{displaymath}
S\subset\Mon(\Lambda_S)=\Mon(\ord(M_1,\dots,M_k))=T(M_1,\dots,M_k)
\end{displaymath}
and the statement \eqref{14.1} is proved.

Statements \eqref{14.2} and \eqref{14.3} follow from proposition
\ref{tv7}, lemma \ref{max} and lemma \ref{rem2}.
\end{proof}
Let $S$ be maximal nilpotent subsemigroup of $\ISA$ of nilpotency
degree $k$ and let $M=M_1\cup M_2\cup\dots\cup M_k$ be the ordered
$A$-partition of $M=\overline{A}^{(1)}\cup A\cup
\overline{A}^{(2)}$, which corresponds to the partial order
$\Lambda_S$. We call the set $(|M_1|, |M_2|,\dots ,
|M_k|)\in\N^{k}$ \textit{ the type of nilpotent subsemigroup  }
$S$ and denote it by $\type(S)$. Set $(|M_1|,\dots
,|M_k|)^{\#}=(|M_k|,\dots ,|M_1|)$.
\begin{proposition} Let $T_1$ and $T_2$ be two maximal nilpotent
subsemigroups of $\ISA$ of nilpotency degree $k$. Then
 \begin{enumerate}[(i)]
\item if $k=2$ then $T_1$ and $T_2$ are isomorphic if and only if
$\type(T_1)=\type(T_2)$ or $\type(T_1)=\type(T_2)^{\#}$.
\item if $k>2$ then $T_1$ and $T_2$ are isomorphic if and only if
$\type(T_1)=\type(T_2)$. $T_1$ and $T_2$ are anti-isomorphic if
and only if $\type(T_1)=\type(T_2)^{\#}$.
\end{enumerate}
\end{proposition}
The proof follows from proposition \ref{mono} and corresponding
statements about nilpotent subsemigroups of $\mathcal{IS}(M)$ from
\cite{GM1} (lemmas 14.1-14.4 and theorem 14.1).
\begin{flushright}
$\square$
\end{flushright}

\vspace{0.1cm} \noindent Department of Mechanics and
Mathematics,\\
Kyiv Taras Shevchenko University,\\
 64, Volodymyrska st., 01033,
Kyiv, UKRAINE,\\ e-mail: {\em gtsyaputa\symbol{64}univ.kiev.ua}

\end{document}